\documentclass[12pt]{article}
\usepackage[english]{babel}
\usepackage{amsmath,amssymb,graphicx,amsthm,enumerate,comment,amscd}
\usepackage[unicode,pdftex]{hyperref}
\usepackage[matrix,arrow,curve]{xy}
\graphicspath{{kartinki/}}

\newtheorem{theorem}{Theorem}

\theoremstyle{remark}

\theoremstyle{definition}

\newcommand{\del}{\smash{\mskip3mu\lower1truept\hbox{\vdots}\mskip3mu}}

\renewcommand{\phi}{\varphi}

\renewcommand{\kappa}{\varkappa}

\newcommand{\sk}{\mathrm{sk}}
\renewcommand{\int}{\mathrm{int}}

\newcommand{\Deg}{\mathrm{Deg}\,}

\begin{document}

\title{Short proof of the Kneser-Edmonds theorem on the degree of a map between closed surfaces}
\author{Andrey Ryabichev%
\footnote{\ IUM; Moscow, Russia
// \texttt{ryabichev@179.ru}
}}
\date{}

\maketitle

\begin{abstract}
Suppose for closed surfaces $M,N$ there exists a continuous map $f:M\to N$ of geometric degree $d>0$.
Then $\chi(M)\le d\cdot\chi(N)$.
This inequality was first proved by Kneser in case of orientable surfaces
and by Edmonds for arbitrary $M,N$.
We give a new simple proof of this result.
Our proof is completely elementary
and does not use additional techniques
(such as the factorisation theorem of Edmonds and the absolute degree theory of Hopf).
\end{abstract}


Let $M,N$ be closed surfaces (i.\,e.\ compact connected $2$-manifolds without boundary, possibly nonorientable).
Given a map $f:M\to N$, the {\it geometric degree $\Deg f$}
is a minimal cardinality of the preimage
of a regular value among all smooth maps homotopic to $f$.
In other words, $\Deg f$ equals the smaller $d\ge0$ such that
there is $f'\sim f$ which is a $d$-sheeted covering over some disk in $N$.

\begin{theorem}\label{th:degree}
Let $f:M\to N$ be a map of geometric degree $d>0$.
Then $\chi(M)\le d\cdot\chi(N)$.
\end{theorem}

In particular, this inequality implies that if $\chi(N)<0$, then a map $M\to N$ is of bounded geometric degree.
Note that if $M,N$ are orientable or, more generally, if $f$ is {\it orientation true},
then its geometric degree equals to usual {\it cohomological degree} (see e.\,g. \cite{brown-schirmer}, \cite{epstein} or \cite{olum}),
but in theorem~\ref{th:degree} the map $f$ need not be orientation true.

This fact was first proved by Kneser \cite{kneser-30} in the case of orientable surfaces.
Edmonds \cite{edmonds} proved it for surfaces with boundary,
considering maps $f$ which restriction to the boundary is $d$-sheeted covering.
His proof in nonorientable case was corrected and improved by Skora \cite{skora}.
The idea is to factorize $f$ up to homotopy as the composition of a pinch map and a branched covering,
see \cite{edmonds} or \cite{skora} for details, see also \cite{ryabichev-degree}.

We give another simple proof of theorem~\ref{th:degree}
which does not use branched coverings and surfaces with boundary
and does not apply the absolute degree theory of Hopf \cite{hopf-28}, \cite{hopf-30}.
We use only the notion of transversality (see e.\,g.\ \cite{hirsch-book})
and some well-known facts on Euler characteristic of surfaces.
The presented proof was obtained by Petr~M.\ Akhmet'ev
and improved by Sergey Melikhov and me.

\begin{proof}
Take a $CW$-structure on $N$ with one $0$-cell $y$ and one $2$-cell.
Clearly, $\sk^1(N)$ is homeomorphic to the bouquet of $k$ circles and $\chi(N)=2-k$.

After a homotopy of $f$ we may assume that it is smooth and that $y$ is a regular value with $d$ preimages.
Moreover, after a small $C^1$-perturbing of $f$
we may assume also  that it is transversal to $\sk^1(N)$.
Then the preimage $f^{-1}(\sk^1(N))\subset M$ is a 
$1$-complex $\Gamma$
and, possibly, a disjoint union of circles $\Omega$.
%
Note that $\chi(\Gamma)=d\cdot(1-k)$ since $\Gamma$ has $d$ vertices
of valence $2k$.
Also note that $y\not\in f(\Omega)$.

We claim that one can homotope $f$ so that it maps the edges of $\Gamma$ to the edges of $\sk^1(N)$ monotonously 
(preserving transversality and fixing $f$ around the vertices of $\Gamma$).
Indeed, otherwise there is an edge $\tilde e\subset\Gamma$ whose image do not cover the corresponding edge $e\subset\sk^1(N)$.
Then we can squeeze $\tilde e$ from $e$ and $y$, and then $y$ will have $2$ preimages less.

Denote by $V\subset N$ a sufficiently small tubular neighborhood of $\sk^1(N)$.
Then the preimage $f^{-1}(V)$ consists of a tubular neighborhood $U_1\supset\Gamma$ and a tubular neighborhood $U_2\supset\Omega$.
Clearly, $U_2$ is a disjoint union of annuli (or empty),
and $U_1$ deformation retracts onto $\Gamma$.
By the claim above, we may assume that
$f|_{\overline{U_1}}:\overline{U_1}\to \overline{V}$ is a $d$-sheeted covering.
Therefore $\partial\overline{U_1}$ has not mote than $d$ components.

The complement $N\setminus V$ is diffeomorphic to a disc, denote it by $D$.
The preimage $f^{-1}(D)$ is a $2$-manifold with boundary, possibly disconnected.
Set $M'=f^{-1}(D)\cup U_2$.
Since $M$ was path-connected,
any component $M'$ has a nonempty boundary.
Therefore, $M'$ has not more than $d$ components, and we have $\chi(M')\le d$.

Finally, note that 
$M'\cup\overline{U_1}=M$
and 
$M'\cap\overline{U_1}$ is a disjoint union of circles.
We obtain
$$\chi(M)=\chi(M')+\chi(\overline{U_1})\le d+d\cdot(1-k)=d\cdot\chi(N),$$
which is what needed to be proven.
\end{proof}

\subsection*{Acknowledgments}
I am very thankful to Petr~M.\ Akhmet'ev and Sergey Melikhov
for their interest for the problem and for sharing the idea of the proof.
I am grateful to Konstantin Shcherbakov for verifying the readability of the text
and for some useful remarks on the proof.
Finally, I want to thank the Independent University of Moscow
and the Young Russian Mathematics award for a financial support.


\end{document}